\documentclass[12pt,final,a4paper]{amsart}
\usepackage{amssymb}
\usepackage{amsmath}
\usepackage{amsfonts}

\newcommand{\Q}{\mathbb{Q}}

\newcommand{\Z}{\mathbb{Z}}

\newcommand{\el}{l}

\newcommand{\lfrac}[2]{(#1)/#2}

\renewcommand{\phi}[0]{\varphi}
\renewcommand{\theta}[0]{\vartheta}

\newcommand{\Span}[1]{\left\langle\, #1 \,\right\rangle}
\newcommand{\Set}[1]{\left\{\, #1 \,\right\}}
\newcommand{\Size}[1]{\lvert\, #1 \,\rvert}

\DeclareMathOperator{\Der}{Der}
\DeclareMathOperator{\Aut}{Aut}
\DeclareMathOperator{\Inn}{Inn}

\DeclareMathOperator{\End}{End}

\newtheorem{dummy}{Dummy}

\numberwithin{dummy}{section}
\numberwithin{equation}{section}

\newtheorem{lemma}[dummy]{Lemma}
\newtheorem{theorem}[dummy]{Theorem}

\theoremstyle{definition}

\theoremstyle{remark}
\newtheorem{rem}[dummy]{Remark}

\begin{document}

\bibliographystyle{amsalpha}

\date{21 December 2005}

\author{A.~Caranti}
\email{caranti@science.unitn.it}

\author{Sandro Mattarei}
\email{mattarei@science.unitn.it}
\address{Dipartimento di Matematica\\
  Universit\`a degli Studi di Trento\\
  via Sommarive 14\\
  I-38050 Povo (Trento)\\
  Italy}

\title[Automorphisms of $p$-groups
  of maximal class]{Automorphisms of $p$-groups\\
  of maximal class}

\begin{abstract}
Juh{\'a}sz has proved that the automorphism group of a group $G$ of
maximal class of order $p^{n}$, with $p \ge 5$ and $n > p + 1$, has
order divisible by  $p^{\lceil(3n-2p+5)/2\rceil}$.

We show that by translating the problem in terms of derivations,
the result can be deduced from the case where $G$ is metabelian.
Here one can use a general
result of Caranti and Scoppola concerning automorphisms of
two-generator, nilpotent metabelian groups.
\end{abstract}



\thanks{Partially supported by MIUR-Italy via PRIN 2003018059
 ``Graded Lie algebras and pro-$p$-groups: representations,
 periodicity and derivations''.}

\maketitle

\thispagestyle{empty}

\section{Introduction}\label{sec:intro}

Baartmans and Woeppel have proved~\cite[Theorem 3.1]{BaaWoe} the
following
\begin{theorem}
  Let $p$ be a  prime, and let $G$ be a $p$-group  of maximal class of
  order $p^{n}$,  which has an abelian maximal  subgroup.  Suppose $G$
  has exponent $p$.  Then $\Aut(G)$ contains a subgroup of order $p^{2
  n - 3}$.
\end{theorem}
The  exponent   restriction  limits  the   size  of  $G$   to  $p^{p}$
(see~\eqref{eq:expo} below).  However, the  main point of  this result
holds    true   more   generally.     Caranti   and    Scoppola   have
proved~\cite{CaSco} (see also~\cite{CaMi}) that any finite, metabelian
$p$-group  has a  subgroup  of  its group  of  automorphisms of  order
$\Size{\gamma_{2}(G)}^{2}$,  where   $\gamma_{2}(G)$  is  the  derived
subgroup of $G$. We thus have in particular
\begin{theorem}\label{thm:CaSco}
  Let $p$ be a prime,  and let $G$ be a $p$-group of maximal
  class of order $p^{n}$, which is metabelian. Then
  $\Aut(G)$ contains a subgroup of order $p^{2 n - 4}$.
\end{theorem}

Juh{\'a}sz has proved in~\cite{Juh} among others the following result.
\begin{theorem}\label{thm:main1}
  Let $p \ge 5$ be a prime, and let $G$ be a $p$-group of maximal
  class of order $p^{n}$, with $n > p + 1$. Then
  $\Aut(G)$ contains a subgroup of order
  $p^{\lceil(3n-2p+5)/2\rceil}$.
\end{theorem}

The aim of this paper is  to show that if one reformulates the problem
in  terms  of derivations, as we do in Section~\ref{sec:derivations},
then  the  general  case of  an  arbitrary
$p$-group of maximal class  of Theorem~\ref{thm:main1} can be shown to
follow  from the  special case  of a  metabelian $p$-group  of maximal
class of Theorem~\ref{thm:CaSco}.

A particular case of Theorem~\ref{thm:main1} has been used by
Malinowska in~\cite{Izabela}.

A  stronger estimate holds  for $3$-groups  of maximal  class, because
these   are   all   metabelian   (see~\cite[Theorem~3.4.3]{L-GMcKay}).
According to  Theorem~\ref{thm:CaSco}, such a group  $G$ of order
$3^n$ has at least $3^{2n-4}$ automorphisms.

After discussing relations between automorphisms and derivations according to our needs in Section~\ref{sec:derivations},
we give our proof of Theorem~\ref{thm:main1} in Section~\ref{sec:maximal}.
The same approach allows us prove also the following result.

\begin{theorem}\label{thm:main2}
  Let $p \ge 5$ be a prime, and let $G$ be a $p$-group of maximal
  class of order $p^{n}$, with $n > p + 1$. Then
  $\Aut(G)$ has an abelian normal subgroup of order $p^{n-2p+7}$.
\end{theorem}

It  must be  noted that  Juh{\'a}sz obtains  his  result in~\cite{Juh}
using  the best  estimate for  the  degree of  commutativity that  was
available  at  the  time  of   his  writing.  We  have  corrected  the
formulation of  Theorem~\ref{thm:main1} to  take account of  the exact
estimate later obtained by Fern{\'a}ndez-Alcober in~\cite{F-A}.

\section{Derivations and automorphisms}\label{sec:derivations}

We conveniently extend to nonabelian groups (written multiplicatively)
a piece of notation usually adopted for the endomorphism ring of
an abelian group (traditionally written additively).
For maps $\sigma,\tau$ from a set $S$ to a (multiplicative) group $G$ we define the map $\sigma+\tau\in G^S$ by setting
$g(\sigma+\tau)=(g\sigma)(g\tau)$ for all $g\in G$.
The ``addition'' operation thus defined, which is not commutative unless $G$ is, makes $G^S$ into a group,
the cartesian product of copies of $G$ indexed by the elements of $S$.
We write as $0$ and $-\sigma$ the identity element and the inverse of $\sigma$ in $G^S$,
and we write $\sigma-\tau$ for $\sigma+(-\tau)$.
The identities $s0=1$ and $s(-\sigma)=(s\sigma)^{-1}$ for $s\in S$ would look more natural
by using the exponential notation $g^\sigma$ for $g\sigma$ which is traditional in similar contexts,
but we avoid doing that to prevent proliferation of exponents.
In the special case where $S=G$, another operation on $G^G$ is given by composition,
written here as (left-to-right) juxtaposition.
It is left-distributive with respect to addition, but not right-distributive, in general.

Recall (cf.~\cite[\S9.5]{L-GMcKay})
that a derivation of a group $G$ into a $G$-bimodule $A$
is a map $\delta: G \to A$
satisfying
$(gh)\delta=(g\delta)h+g(h\delta)$ for all $g,h\in G$.
The set of derivations of $G$ into $A$, denoted by $\Der(G,A)$,
is an abelian group with operation induced by that on the codomain $A$, as described above.
We define the {\em kernel} of a derivation $\delta$ as $\ker\delta=\{g\in G: g\delta=0\}$,
bearing in mind that this is a subgroup of $G$ but need not be normal.

\begin{rem}\label{aec:002}
If  $N$ is  a normal  subgroup of  $G$ for  which $A$  is  the trivial
bimodule then we  can and will identify $\Der(G/N,A)$  with the subset
of  $\Der(G,A)$ consisting  of the  derivations whose  kernel contains
$N$.
\end{rem}

Now let $A$ be an abelian normal subgroup of a group $G$.
Then it is customary to make $A$ into a $G$-bimodule with the trivial action on the left
and the conjugation action on the right, that is,
$g\cdot a\cdot h=a^h$ for $a\in A$ and $g,h\in G$.
In this context, writing the group operation in $A$ multiplicatively, as in $G$,
the condition for $\delta:G\to A$ being a derivation reads
$(gh)\delta=(g\delta)^h(h\delta)$.
In particular, this readily implies that
$1\delta=1$ and $g^{-1}\delta=((g\delta)^{-1})^{g^{-1}}$.
It follows that if $\delta\in\Der(G,A)$ then the map $\alpha=1+\delta$, given by
$g\alpha=g(g\delta)$ according to notation introduced earlier, is an endomorphism of $G$.
Conversely, if $\alpha$ is an endomorphism of $G$ which sends $A$ into
itself  and induces  the  identity  map on  the  quotient $G/A$,  then
$-1+\alpha$        belongs        to        $\Der(G,A)$.         Since
$(1+\delta_1)(1+\delta_2)=1+\delta_1+\delta_2+\delta_1\delta_2$     for
$\delta_1,\delta_2\in\Der(G,A)$,   the   operation  ``$\bullet$''   on
derivations                          given                          by
$\delta_1\bullet\delta_2=\delta_1+\delta_2+\delta_1\delta_2$      makes
$\Der(G,A)$  into  a  monoid,  and the  correspondence  $\delta\mapsto
1+\delta$  becomes  a  monomorphism  into the  monoid  $\End(G)$  with
respect to composition.  We record part of these conclusions for later
reference.

\begin{lemma}\label{lemma:monoids}
Let $A$ be an abelian normal subgroup of a group $G$.
The map sending $\delta$ to $1+\delta$
is a monomorphism of the monoid
$\Der(G,A)$ with the operation $\bullet$
into the monoid $\End(G)$ with respect to composition.
Its image consists of the endomorphisms of $G$ which send $A$ into itself
and induce the identity map on $G/A$.
\end{lemma}

The endomorphism $1+\delta$ is injective provided $\delta$ maps no element $g$ of $G$ (or, equivalently, of $A$) to its inverse.
A sufficient condition for this to occur is, for example, that $\delta$ is nilpotent,
in the sense that some power $\delta^k$ vanishes,
because then $1+\delta$ admits the inverse $1-\delta+\delta^2-\delta^3+\cdots$.
(This occurs in Lemma~\ref{lemma:groups} below, with $k=2$.)
In case $G$ is a finite group,
a sufficient condition for $1+\delta$ being an automorphism of $G$ for all $\delta\in\Der(G,A)$,
and hence for $1+\Der(G,A)$ to be a subgroup of $\Aut(G)$,
is that $A$ is contained in the Frattini subgroup of $G$.
In fact, in that case the image of $1+\delta$ supplements the Frattini subgroup,
whence $1+\delta$ is surjective, and thus injective by finiteness of $G$.

On the subset $\Der(G/A,A)$ of $\Der(G,A)$ the operation $\bullet$ coincides with addition,
because $\delta_1\delta_2=0$ for $\delta_1,\delta_2\in\Der(G/A,A)$.
In particular, it is commutative in this case.
The properties of the correspondence $\delta\mapsto 1+\delta$ stated in Lemma~\ref{lemma:monoids}
read as follows.

\begin{lemma}\label{lemma:groups}
Let $A$ be an abelian normal subgroup of a group $G$.
The map sending $\delta$ to $1+\delta$
is a monomorphism of the additive group
$\Der(G/A,A)$ into $\Aut(G)$.
Its image is the abelian subgroup consisting of the automorphisms which send $A$ into itself and induce the identity
map on $G/A$.
\end{lemma}

A familiar instance of Lemma~\ref{lemma:groups} is when $A$ is the centre of $G$.
Then derivations $\delta\in\Der(G/A,A)$ are the same thing as group homomorphisms of $G/A$ into $A$,
and the correspondence $\delta\mapsto 1+\delta$ maps $\Der(G/A,A)$ onto the group of
{\em central} automorphisms of $G$.

We will need the following fact on derivations of a group $G$ into an (abelian) term of its lower central series
$G=\gamma_1(G)\ge\gamma_2(G)\ge\cdots$.

\begin{lemma}\label{lemma:down}
  Suppose that $\gamma_r(G)$ is abelian and let $\delta\in\Der(G,\gamma_r(G))$.
  Then $\gamma_i(G)\delta\subseteq \gamma_{i+r-1}(G)$ for all $i\ge 1$.

  In particular if $G$ is nilpotent, with $\gamma_{n}(G) = 1$, then we
  have, according to Remark~\ref{aec:002},
  \begin{equation*}
    \Der(G, \gamma_{r}(G))
    =
    \Der(G / \gamma_{n - r + 1}(G), \gamma_{r}(G)).
  \end{equation*}
\end{lemma}

\begin{proof}
Since $1+\delta$ is an endomorphism of $G$, for all $g,h\in G$ we have
\[
[g,h]([g,h]\delta)=[g,h](1+\delta)=
[g(1+\delta),h(1+\delta)]=[g(g\delta),h(h\delta)],
\]
and hence
$[g,h]\delta=[g,h]^{-1}[g(g\delta),h(h\delta)]$.
The commutator identity
\[
[gu,hv]=[g,v]^u[g,h]^{vu}[u,v][u,h]^v
\]
shows that
\[
[g,h]^{-1}[gu,hv]\in\gamma_{i+j+r-1}(G)
\]
if $g\in\gamma_i(G)$, $u\in\gamma_{i+r-1}(G)$, $h\in\gamma_j(G)$ and $v\in\gamma_{j+r-1}(G)$.
Since $\gamma_{i+1}(G)$ is generated by all commutators $[g,h]$ with
$g\in\gamma_i(G)$ and $h\in G=\gamma_1(G)$,
our claim follows by induction on $i$ by taking $u=g\delta$ and $v=h\delta$.
\end{proof}

\section{Automorphisms of $p$-groups of maximal class}\label{sec:maximal}

We take~\cite{L-GMcKay} as a  reference for $p$-groups of
maximal class, but see also~\cite[III.14]{Hup}, and  Blackburn's
original paper~\cite{black}.

In this  section, $G$ will  be a $p$-group  of maximal class  of order
$p^n$, with $p \ge 5$ and $n \ge 4$.
As usual, write $G_i=\gamma_i(G)$ for $i\ge 2$, and define a maximal
subgroup $G_{1}$ of $G$ by
\begin{equation*}
  G_{1}
  =
  C_{G}(G_{2} / G_{4})
  =
  \{g\in G\mid [G_2,g]\subseteq G_4\}.
\end{equation*}
In particular, $G_{n-1} > G_n = 1$.

The degree of  commutativity $l$ of $G$ is defined as $n-3$
if $G_{1}$
is abelian, and otherwise as the largest integer
$l$ such that $[G_{i}, G_{j}] \le G_{i+j+l}$  for all $i, j \ge 1$.
Since $[G_1,G_1]=[G_1,G_2]$ we have $l\le n-3$ in all cases.
One can show (\cite[Theorem~3.3.5]{L-GMcKay},
\cite[Hauptsatz~III.14.6]{Hup}) that for $n > p + 1$ the degree of
commutativity of a group
$G$ of maximal class of order $p^n$ is  positive, that is,
\begin{equation*}
  G_{1}
  =
  C_{G}(G_{i} / G_{i+2})
\end{equation*}
for all $i = 2, \dots , n-2$.
From now on we take $n > p + 1$.
Choose $s_{1} \in G_{1} \setminus G_{2}$
and $s \in  G \setminus G_{1}$,
and define $s_{i+1} =  [s_{i}, s]$ for $i\ge 1$.
We then have $G_{i}  = \Span{s_{i}, G_{i+1}}$ for $i = 1,\ldots, n-1$.

\begin{lemma}\label{lemma:kernel}
  Let $r \ge (n-l)/2$ and $\delta \in \Der(G,G_r)$.
  Then
  \begin{equation*}
    G_{n-r+1} \le \ker(\delta),
  \end{equation*}
  and hence $\delta$ can be viewed as
  an element of
  \begin{equation*}
    \Der(G/G_{n-r+1},G_r).
  \end{equation*}
\end{lemma}

\begin{proof}
  Note that $G_{r}$ is abelian, as $[G_{r}, G_{r}] \le G_{2 r + l} \le
  G_{n} = 1$.
  The conclusion follows at once from  Lemma~\ref{lemma:down}.
\end{proof}

Let now $G'$ be a group of maximal class which is metabelian, that is,
with the obvious notation, $[G_{2}', G_{2}'] = 1$. A result
of~\cite{CaSco} (see also~\cite{CaMi}) guarantees that $G'$ has
plenty of automorphisms.
\begin{theorem}
  \label{thm:metabelian}
  Let $M = \Span{x, y}$ be a metabelian, $2$-generator finite
  nilpotent group, and let
  $M_{2} = [M, M]$ be its derived subgroup.
  Then for all $u, v  \in M_{2}$ there is an automorphism of
  $M$ such that
  \begin{equation*}
    \left\{
    \begin{aligned}
      x &\mapsto x \cdot u,\\
      y &\mapsto y \cdot v.
    \end{aligned}
    \right.
  \end{equation*}
\end{theorem}

With the terminology introduced in the previous section we can rephrase the conclusion of Theorem~\ref{thm:metabelian} as follows:
for all $u,v\in M_{2}$ there is a derivation $\delta\in\Der(M,M_2)$ such that
$x\delta=u$ and $y\delta=v$.

We intend to exploit this in the following way. Given an arbitrary
$p$-group $G$  of maximal class  of order
$p^n$, with $n > p + 1$, we will show that there are a suitable $r$,
and a suitable metabelian $p$-group $G'$ of maximal class, of the same
order as $G$, such that
\begin{itemize}
\item $G_{r}$ is abelian,
\item $G / G_{n-r+1}$  is isomorphic to $G' / G_{n-r+1}'$, and
\item the $G / G_{n-r+1}$-module $G_{r}$ is similar to the $G' /
  G_{n-r+1}'$-module $G_{r}'$.
\end{itemize}
It will follow that
\begin{equation}\label{aec:006}
  \Der(G / G_{n-r+1}, G_{r})
  \cong
  \Der(G' / G_{n-r+1}', G_{r}').
\end{equation}
Now Theorem~\ref{thm:metabelian} tells us that $G'$ has many automorphisms,
that is, the set at the right hand side of~\eqref{aec:006} is large, so that the
set at the left hand side is also large, and $G$ in turn has many automorphisms.

We begin by defining $G'$, following~\cite[p.~83--84]{black}, by the
presentation
\begin{align*}
  G' = \langle
  s', s_{i}', i = 1, \dots, n - 1:
  &\
  s'^{p} = 1,
  \\&
    [s_{i}', s'] = s_{i+1}'\ \text{for $i = 1, \dots, n-2$},
    \\&
      [s_{i}', s_{j}'] = 1\ \text{for $i, j = 1, \dots, n-1$},
      \\&
      {s_{i}'}^{p} {s_{i+1}'}^{\binom{p}{2}} \cdots s_{i+p-1}' = 1
      \ \text{for $i = 1, \dots, n$}
      \rangle.
\end{align*}
(We assume $s'_{i} = 1$ for $i \ge n$.)
This group may be constructed in the following way. One starts with
the abelian group
\begin{align*}
  M = \langle
  s_{i}', i = 1, \dots, n-1 :
  &\
  [s_{i}', s_{j}'] = 1\ \text{for $i, j = 1, \dots, n-1$},
  \\&
  {s_{i}'}^{p} {s_{i+1}'}^{\binom{p}{2}} \cdots s_{i+p-1}' = 1
  \ \text{for $i = 1, \dots, n-1$}
  \rangle.
\end{align*}
This has order $p^{n-1}$, and its structure can be understood by
reading the last group of relations backwards. Now $M$ admits
an automorphism $\sigma : s'_{i} \mapsto s'_{i} s'_{i+1}$, as
$\sigma$ preserves the defining relations. Moreover, for $i \ge 2$ one
has, by~\cite[Hilfssatz~10.9(b)]{Hup} or~\cite[Corollary~1.1.7]{L-GMcKay},
\begin{equation*}
  \sigma^{p}(s'_{i-1})
  =
  s'_{i-1}
  \cdot
  {s_{i}'}^{p} {s_{i+1}'}^{\binom{p}{2}} \dots s_{i+p-1}'
  =
  s'_{i-1},
\end{equation*}
so that $\sigma$ has order $p$.
Now $G'$ above can be constructed as
the cyclic extension of $M$ by a cyclic group $\Span{s'}$ of order $p$,
where $s'$ induces $\sigma$ on $M$.

Take $r= n - \el - 1$, and $A = G_{r}$.
We have $[G_{1}, A] = 1$, so that in
particular $A$ is abelian.
Note that $r>(n-l)/2$ because $l\le n-3$, and hence
$\Der(G,A)=\Der(G/G_{\el+2},A)$ according to Lemma~\ref{lemma:kernel}.

It is time to recall some basic facts about $p$-groups of maximal class,
valid under our assumption $n>p+1$.
If $g \in G \setminus G_{1}$, then $g \notin
C_{G}(G_{i}/G_{i+2})$, for $i = 1, \dots, n-2$. Thus $C_{G}(g) \cap
G_{1} = G_{n-1} = Z(G)$.
It follows that $C_{G}(g) = \Span{g,
  G_{n-1}}$, so that $g^{p} \in G_{n-1}$. Also, the conjugacy class
$g^{G}$ of $g$ has order $\Size{G : C_{G}(g)} = p^{n-2}$, so that
$g^{G} = g G_{2}$.

As $s  \notin G_{1}$, we obtain  in particular that $s$  and $s s_{i}$
are conjugate, for  $i \ge r$, and hence the  elements $s^{p}$ and
$(ss_{i})^{p}$,  which  by  the  above  lie  in the  centre  of  $G$,  do
coincide. If $i \ge r$, we  have that $s_{i}$ commutes with all of the elements
$s_{j} = [s_{i},  \underbrace{s, \dots, s}_{j-i}]$, for $j  \ge i$.
Consequently, we have
\begin{equation}\label{eq:expo}
  1
  =
  s^{-p} (s s_{i})^{p}
  =
  s_{i}^{p} s_{i+1}^{\binom{p}{2}} \dots s_{i+p-1},
\end{equation}
again by~\cite[Hilfssatz~10.9(b)]{Hup} or~\cite[Corollary~1.1.7]{L-GMcKay}.
These relations define $G_{r}$ as an abelian group generated by
the $s_{i}$, for $i \ge r$, so that $G_{r}$ is isomorphic to $G'_{r}$.

Because $[G_{1}, G_{1}] \le G_{\el + 2} = N$, the quotient $G_{1}/N$ is abelian.
As above, we have
\begin{equation}
  \label{aec:007}
  s_{i}^{p} s_{i+1}^{\binom{p}{2}} \dots s_{i+p-1}
  \equiv
  1 \pmod{N}
\end{equation}
for all $i \ge 2$. Since $s^{p} \equiv (s s_{1})^{p} \equiv 1 \pmod{N}$,
equation~\eqref{aec:007} also holds for $i = 1$. We thus see that $G
/ N = G / G_{\el + 2}$ is isomorphic to the corresponding factor group
$G' / G'_{\el + 2}$ of $G'$.
Finally, the action of $G$ on $G_{r}$ is given by $[s_{i}, s_{1}]= 1$
and $[s_{i}, s] = s_{i+1}$, for $i \ge r$.
Thus, the $G / G_{\el +2}$-module $G_{r}$
is similar to the $G' / G'_{\el +2}$-module $G'_{r}$.

According to Theorem~\ref{thm:metabelian} and Lemma~\ref{lemma:monoids},
we conclude that for all $u,v\in A$ there is an automorphism of $G$ determined by
\begin{equation*}
  \phi_{u,v} :
  \left\{
  \begin{aligned}
    s &\mapsto s \cdot u,\\
    s_{1} &\mapsto s_{1} \cdot v.
  \end{aligned}
  \right.
\end{equation*}
A comment following Lemma~\ref{lemma:monoids} shows that these automorphisms form a subgroup of $\Aut(G)$.
In particular, the automorphisms among these which fix $s$ form a subgroup
$\mathcal{H}=\Set{\phi_{1, v} : v \in A}$, of order $p^{n-r}$.

Suppose that $\phi_{1,v} \in \Inn(G)$ for  some $v\neq 1$, so that $s^{g} = s$ and
$s_{1}^{g} =  s_{1} v \ne s_{1}$ for some $g \in G$.
Since $C_{G}(s)= \Span{s, G_{n-1}}$, we have that  $g \equiv s^{i}$ modulo the centre $G_{n-1}$
for some $i$, and $v = [s_{1}, s^{i}] \in G_{2} \setminus G_{3}$.
Thus the group $\mathcal{H}$ intersects $\Inn(G)$ trivially if $r  > 2$.
We will verify below that the metabelian case $r = 2$, that is $l = n -3$,
is covered by Theorem~\ref{thm:metabelian}, so in the following we assume
$l < n - 3$.

It follows that $\Aut(G)$ contains the subgroup $\mathcal{H} \Inn(G)$,
with $\mathcal{H} \cap \Inn(G) = \Set{1}$.
Note that all automorphisms $\phi_{u,v}$ belong to $\mathcal{H} \Inn(G)$.
In fact, for $u \in A$ one
has that $s u \in s A
\subseteq s G_{2}$, so $s u$ is conjugate to $s$. If $s u = s^{g}$ for
some $g \in G$, then composing $\phi_{u,v}$ with conjugation by
$g^{-1}$ one obtains an element of $\mathcal{H}$.

So we get that $\Aut(G)$ has a subgroup of size at least $p^{c}$, with
\begin{equation*}
  c \ge n - 1 + n - r =  n + \el.
\end{equation*}
Using the estimate $2l\ge n-2p+5$ of~\cite{F-A}
(see~\cite[Theorem~3.4.11]{L-GMcKay} for a version of this bound weakened
by one)
we conclude that
\begin{equation*}
  c \ge \dfrac{3 n - 2 p + 5}{2},
\end{equation*}
thus completing a proof of Theorem~\ref{thm:main1} except for the case where $G$ is metabelian.
However, in that case Theorem~\ref{thm:metabelian} provides $p^{2(n - 2)}$ distinct automorphisms of $G$,
and this number exceeds $p^{\lfrac{3 n - 2 p +5}{2}}$ for $p \ge 3$.

\begin{rem}
The spirit of the times suggests an alternative description of $G'$,
as in~\cite[Examples~3.1.5]{L-GMcKay}, which we sketch here.
Let $\theta$ be a primitive $p$th root of unity over the rational field $\Q$.
Then the abelian group $M$ can be realised (in additive notation) as the additive group of the quotient ring
$\Z[\theta]/(\theta-1)^{n-1}$,
where the residue class of $(\theta-1)^{i-1}$ plays the role of $s'_i$.
(The defining relations in Blackburn's presentation for $M$ are then all consequences of the relation
$(1+(\theta-1))^p-1=0$.)
We construct $G'$ as the cyclic extension of $M$ by a cyclic group $\Span{s'}$ of order $p$,
where $s'$ acts on $M$ by multiplication by $\theta$.
Now, the derivations $\delta\in\Der(G',M)$ such that $s'\delta=0$
correspond to the endomorphisms of $M$ as $\langle  s'\rangle$-module,
which are clearly given by all multiplications by polynomials in $\theta$.
In particular, this gives an explicit description of $\Der(G',A)$
which allows one to construct $\mathcal{H}$ without recourse to Theorem~\ref{thm:metabelian}.
\end{rem}

In order to prove Theorem~\ref{thm:main2}, consider the subgroup
$G_t$ of $G_r$, where $t=\max(n-l-1,\lceil\frac{n+1}{2}\rceil)$.
According to Lemma~\ref{lemma:kernel}, and because $n-t+1\le t$,
we have
\[
\Der(G,G_t)=\Der(G/G_{n-t+1},G_t)
=\Der(G/G_t,G_t).
\]
Lemma~\ref{lemma:groups} implies that
$\Set{\phi_{u, v} : u,v \in G_t}$
is an abelian subgroup of $\Aut(G)$ isomorphic with the additive group $\Der(G/G_t,G_t)$,
and hence of order
$|G_t|^2=p^{2(n-t)}$.
This abelian subgroup of $\Aut(G)$ is normal, again according to Lemma~\ref{lemma:groups},
because $G_t$ is a characteristic subgroup of $G$.
Because $2l\ge n-2p+5$ and $p\ge 5$ we see that
$2(n-t)\ge n-2p+7$, thus proving Theorem~\ref{thm:main2}.

\bibliography{References}

\end{document}